\documentstyle[11pt,amscd,amsfonts]{amsart}
\newtheorem{thm}{Theorem}[section]
\newtheorem{lem}[thm]{Lemma}
\newtheorem{cor}[thm]{Corollary}
\newtheorem{prop}[thm]{Proposition}
\newtheorem{notitle}[thm]{ }

\newtheorem{problem}[thm]{Problem}

\theoremstyle{definition}

\def \smash {\wedge}

\def \L {\mathop{\mathrm {L}}}

\def \C {\mathbb{C}}
\def \F {\mathbb{F}}
\def \R {\mathbb{R}}
\def \T {\mathbb{T}}
\def \Z {\mathbb{Z}}

\def \CA {{\cal A}}

\def \CC {{\cal C}}
\def \CG {{\cal G}}

\def \sc {spin^c}

\def \a {\alpha}

\def \d {\delta}
\def \f {\phi}

\def \lam {\lambda}

\def \s {\sigma}

\begin{document}

\baselineskip.525cm

\title[Cohomotopy invariants: II]{A stable cohomotopy refinement of\\ Seiberg-Witten invariants: II}
\author[Stefan Basuer]{Stefan Bauer} \address{Fakult\"at f\"ur Mathematik,
Universit\"at Bielefeld \newline \hspace*{.375in}PF 100131, D-33501 Bielefeld}
\email{bauer@@mathematik.uni-bielefeld.de}
\begin{abstract} A gluing theorem for
the stable cohomotopy invariant defined in the 
first article in this series of two
gives new results on diffeomorphism types of 
decomposable manifolds.
\end{abstract}
\maketitle


\section{Introduction\label{Intro}}

A stable cohomotopy refinement of Seiberg-Witten invariants was defined in
the first article \cite{BauerFuruta} of this series of two. New invariants 
have to stand some acid tests: Do they give new insights?  Can they be computed?
The present article addresses these questions. It turns out that the
stable cohomotopy invariants actually encode more information on 
four-manifolds 
than the integer valued Seiberg-Witten invariants. 
Recall from \cite{BauerFuruta} that for a $K$-oriented, closed Riemannian
4-manifold the monopole map $\mu: \CA\to \CC$ is an ${\T}$-equivariant
fiber preserving map between Hilbert space bundles over
$Pic^0(X)=H^1(X;{\R})/H^1(X;\Z)$, which defines an element $[\mu]$ in an
equivariant stable cohomotopy group  $\pi^b_{{\T},H}(Pic^0(X);{\mathrm
ind}(D))$. Here $b=b^+_2(X)$ denotes the maximal dimension of a positive definite
subspace of second de Rham cohomology of $X$ and $ind (D)$ denotes the virtual index
bundle over $Pic^0(X)$ of the twisted Dirac operator. The subscript $H$
denotes a universe for the action of the group ${\T}$ of complex numbers of
unit length. This technical device is
needed in the equivariant stable homotopy setting to specify the group action
on the suspension coordinates. The key result of this
article is a connected--sum theorem for this invariant.  Recall that a
$K$-theory orientation of a connected sum of 4-manifolds  uniquely induces
$K$-theory orientations on the respective summands.

\begin{thm} \label{Summensatz}
For a connected sum $X=X_0\#X_1$ of 4-manifolds, the stable
            equivariant cohomotopy invariant
            is the smash product of the 
            invariants of its summands 
            \[[\mu_X]=[\mu_{X_0}]\smash[\mu_{X_1}].\]
        \end{thm}

Stating it loosely, the monopole map $\mu_X$ has the same 
stable cohomotopy invariant as the product $\mu_{X_0}\times \mu_{X_1}$.

Let's consider some applications of this theorem. In these applications rudimentary facts 
about equivariant stable cohomotopy groups already gain geometrically significant old and
new information about four-manifolds.

First some implications which are well known from Seiberg-Witten theory:
If  both summands have nonzero $b^+_2$, then by dimension reasons the smash product of the 
cohomotopy invariants lies in the kernel of the homomorphism to the integers
comparing the stable cohomotopy invariants to the Seiberg-Witten invariants.
Thus one regains the folklore theorem stating that the
Seiberg-Witten invariant for a connected sum of $4$-manifolds
vanishes if both summands have nonzero $b^+_2$.
The blowup formula for Seiberg-Witten invariants
follows from the above theorem by applying a well known result about
the degree of ${\T}$-equivariant self maps of representation spheres.

If the manifold
$X$ has vanishing first Betti number, then by forgetting the equivariance, the
stable cohomotopy invariant gives rise to an element of the stable $k$-stem $\pi^s_k$ with 
$k= {\mathrm ind}_{\R}(D)-b^+_2$. The above theorem then states  that
the invariant is multiplicative under connected sums. 
If the $\sc$-structure on $X$ comes from an almost complex structure, 
then the invariant is an element in the group $\pi^s_1\cong {\Z}/2$, with the
Hopf map representing the nontrivial element. It turns out that the invariant
is the Hopf map if and only if both the integer 
Seiberg-Witten invariant is odd and $b^+_2\equiv 3\, mod\, 4$. The fact that 
the cube of the Hopf map is nonvanishing, together with known computations
of the Seiberg-Witten invariants, lead to the following sample applications:

\begin{cor}\label{Separatisten2}
            The connected sum of two symplectic four-manifolds with 
            vanishing first Betti number and with $b^+_2\equiv
            3\,mod \,4$ each 
            does not split off a manifold
            with $b^+_2\equiv 1\, mod\, 4$.
          \end{cor}

\begin{cor}\label{Separatisten3}
            Let $X$ be a connected sum of three symplectic four-manifolds with 
            vanishing first Betti number and with $b^+_2\equiv
            3\,mod \,4$ for each summand.  
            Suppose $X$ splits as a connected sum $X\cong X_1\# X_2$ 
            with $b^+_2(X_1)\equiv 1\, mod\, 4$. Then the
            intersection form on $X_2$ is negative definite.
          \end{cor}

\begin{cor} \label{K3}
            Let $K$ denote the K3-surface and suppose there is an oriented
            diffeomorphism $f:X_1\#K\#K\to X_2\#K\#K$, 
            where the $X_i$ are simply connected K\"ahler manifolds
            with $b^+_2= 3\, mod\, 4$. 
            Then the integer Seiberg-Witten invariants
            of $X_1$ and $X_2$ are the same $mod\, 2$. More precisely,
            let the basic set $B(X_i)\subset H^2(X_i;{\Z})$ consist of the characteristic classes 
            with odd Seiberg-Witten number. If one views $H^2(X_i;{\Z})$ as a
            direct summand in $H^2(X_i\#K\#K;{\Z})$, then $f$ preserves the basic sets,
            $f^*B(X_2)=B(X_1)$. 
        \end{cor}

\begin{cor} \label{elliptic}
            Suppose the connected sum $\#_{i=1}^mE_i$ of simply connected
            minimal elliptic surfaces of odd 
            geometric genus is diffeomorphic to a 
            connected sum $\#_{j=1}^nF_j$ of elliptic surfaces. If $m < 4$, 
            then $n=m$ and the $F_j$ and the $E_i$ are 
            diffeomorphic up to permutation.
        \end{cor}

By the computation in Lemma 3.5 of \cite{BauerFuruta}, the
connected sum of five almost complex manifolds with vanishing
first Betti number always has vanishing stable cohomotopy invariant.
Mikio Furuta suggested that in the case of four summands
there still might exist nonvanishing invariants. Indeed,
one has the following result, for which
the equivariant setting is used in an essential way: 

\begin{cor} \label{power4}
Let $X=\#_{i=1}^4X_i$ be the connected sum of four almost complex manifolds with
vanishing first Betti numbers. Then the stable cohomotopy invariant is nonvanishing
if and only if all of the following conditions hold: For every summand $X_i$,
the integer Seiberg-Witten invariant is odd and $b^+_2\equiv 3\, mod\, 4$. For $X$,
furthermore, the congruence $b^+_2\equiv 4\, mod \, 8$ holds. \\In particular,
the statements in 1.4 and 1.5 also hold for 
the connected sum of four manifolds, as long as the resulting manifold
satisfies the congruence $b^+_2\equiv 4\, mod \, 8$.
\end{cor}

\noindent
{\it Acknowledgement:}
I am grateful to Simon Donaldson for a preliminary version of 
\cite{Donaldson} and to an unnamed referee, who
pointed out the use of an incorrect lemma from the literature 
at a central point in an earlier version of this paper \cite{Bauer}. 
Conversations with Kim Froyshov and Peter Kronheimer were 
helpful in filling the gap and smoothing the argument. 

\section{The setup}

Let $X$ be the disjoint union of a finite number, say $n$, of closed connected 
Riemannian $4$-manifolds $X_i$, each equipped with a $K$-theory orientation. Suppose each component
contains a separating ``long neck'' $N(L)_i=[-L,L]\times S^3$. So it is
a union  \[X_i=X_i^-\cup X_i^+\] of closed submanifolds with
common boundary $\partial X_i^\pm=\{0\}\times S^3$. 
The length $2L\gg 8$ of the long neck will be specified below, its radius is assumed to be
equal in all components.

For an even permutation $\tau\in A_n$, let $X^\tau$ be the manifold obtained
from $X$ by interchanging the positive parts of its components, that is
\[{X}_i^\tau=X_{i}^-\cup X_{\tau(i)}^+.\]
The $K$-theory orientation of $X$ induces by gluing a $K$-theory orientation of $X^\tau$.
The main result of the present article compares the stable cohomotopy invariants
of the two manifolds $X$ and $X^\tau$. The connected sum theorem will be an immediate consequence.

For the manifold $X$ the monopole map of \cite{BauerFuruta} can be described as follows: 
Let $S^+$ and $S^-$  denote the Hermitian 
rank-2 bundles associated to the given $K$-orientation and let $A$ denote a 
$\sc$-connection, which we may assume to induce the flat connection on $det(S^\pm)$ over
the long neck. Fix once and for all identifications
of the spinor bundles and the chosen $\sc$-connections over the $n$ copies of 
$[-L,L]\times S^3$ in $X$.

The gauge group ${\cal G}=map(X,{\T})$
acts on the space $\Gamma(S^\pm)$ of spinors via multiplication with $u:X\to {\T}$, 
on connections via addition of $i\ u{\mathrm d}u^{-1}$ and trivially on forms. 
Consider the subgroup ${\CG}_r\subset {\CG}$ consisting of gauge transformations
which are trivial over the ``short neck'' $N(1)=[-1,1]\times \coprod_{i=1}^n S^3$.
The group ${\CG}_r$ decomposes into a product
of gauge groups, each corresponding to one of the manifolds $X_i^\pm$.
The orbit of the action of ${\CG}_r$ on the space of $\sc$-connections is of the 
form $(A+i\,{\mathrm ker}\,{\mathrm d}_r)$, where
${\mathrm ker}\,{\mathrm d}_r\subset {\mathrm ker}\, {\mathrm d}$ is the space of closed 
1-forms on $X$ vanishing identically on the short neck.  Using the
identification of the chosen $\sc$-connections $A$ and $A^\tau$ over the short neck,
the space $A+i\, {\mathrm ker}\, {\mathrm d}_r \cong A^\tau +i\, {\mathrm ker}\, {\mathrm d}_r$ 
can be viewed as a subspace of the space of $\sc$-connections both over $X$ and $X^ \tau$.
After suitable Sobolev completion, $A+i\, {\mathrm ker}\, {\mathrm d}_r/{\CG}_r$ identifies this
way both with $Pic^0(X)=H^1(X;{\R})/H^1(X;{\Z})$ and $Pic^0(X^\tau)$.

Let $\CA$ and $\CC$ denote the quotients
\begin{eqnarray*}
 {\CA}&=&(A+i\, {\mathrm ker}\, {\mathrm d}_r)\times\left(\Gamma(S^+)\oplus \Omega^1(X)\right)/{\cal G}_r\\
{\CC}&=&(A+ i\, {\mathrm ker}\, {\mathrm d}_r)\times
                \left(\Gamma(S^-)\oplus \Omega^2_+(X)\oplus H^1(X;{\R})\oplus
                \Omega^0(X)/{\R}\right)/{\cal G}_r
\end{eqnarray*}
by the action of the gauge group. Both spaces are bundles over $Pic^0(X)$ and the monopole map
\[\mu=\widetilde{\mu}/{\cal G}_0:{\CA}\to{\CC}\]
is a fiber preserving, ${\T}$-equivariant map over $Pic^0(X)$, defined by
\[  (A',\phi, a)\mapsto (A',\,D_{A'+a}\phi,\,F_{A'+a}^+-\s(\phi),pr(a),
            \,{\mathrm d}^*a).
\]
Here $\s(\phi)$ denotes the trace free endomorphism $\phi\otimes\phi^*-
{1\over 2}\vert\phi\vert^2\cdot \mathop{\mathrm{id}}$ of $S^+$, considered as 
a selfdual 2-form on $X$. The map $pr$ is a linear map $\Omega^1(X)\to H^1(X;{\R})$,
which is an isomorphism on the space ${\mathrm ker}({\mathrm d}^*+{\mathrm d}^+)$ of harmonic one-forms. 
Any such map will do. For the purposes of this paper, we make a choice as follows: 
Fix $b_1$ closed curves, smoothly embedded in the complement
of the long neck, which form a basis of $H_1(X;{\R})$. For a continuous one-form $a$ on $X$, integration
of $a$ along multiples of these cycles defines the element $pr(a)\in Hom(H_1(X;{\R}), {\R})=H^1(X;{\R})$.

It is convenient to use fiberwise $L^2_k$ and $L^2_{k-1}$ Sobolev completions ($k\ge 4$)
for $\CA$ and $\CC$, respectively.
The fiberwise ${\T}$-action on $\CA$ and $\CC$ is trivial on forms and is 
given by complex multiplication on spinors.

To compare the monopole maps $\mu: \CA\to\CC$ on $X$ and 
$\mu^\tau:\CA^\tau\to\CC^\tau$ on $X^\tau$, we will use gluing maps
$\CA\to \CA^\tau$ and $\CC\to\CC^\tau$.
To define these maps, we first need a smooth path
\[\psi:[0,1] \to SO(n)\]
starting from the unit, i. e. $\psi(0)= \mathrm{id}$, and ending 
at $\tau$, considered as the permutation matrix 
$(\d_{i, \tau(j)})_{i,j}\in SO(n)$. 
A second ingredient in the construction is a smooth function \[\gamma:[-L,L]\times S^3\to [0,1],\]
depending only on the first variable. This function $\gamma$ is supposed to vanish
on the $[-L, -1]$-part of the neck and to be identical to $1$ on the $[1,L]$-part. 

For a section $e$ of a bundle $E$ over $X$, denote by $e_i$ its restriction to the bundle
$E_i= E_{\vert X_i}$. Suppose 
the restrictions of $E_i$ to the long necks are identified with a bundle
$F$. Using these identifications, the restrictions $E_{\vert X_i^\pm}$ glue together
to give a bundle $E^\tau$ over $X^\tau$.
Smooth sections of $E$ will be patched together to give smooth sections of $E^\tau$
as follows: The restrictions of $e_i$ to the complement of the long neck 
$N(L)$ remain unchanged.  Over the long neck, the
restrictions of the sections $e_i$ can be viewed as the components of a section
\[
\Vec{e} = \left(\begin{matrix}e_1\cr\vdots\cr e_n\end{matrix}\right)
\]
of the bundle $\oplus_ {i=1}^n F$ over $[-L,L]\times S^3$. The $i$-th component of the section
\[
\Vec{e^\tau} = (\psi\circ\gamma) \cdot\Vec{e}
\]
now restricts to the section $e_i$ over the $[-L,-1]$-part
and to the section $e_{\tau(i)}$ over the $[1, L]$-part. Patching 
parts together, we obtain a smooth section $e^\tau$ of the bundle $E^\tau$.

This gluing construction, applied to forms $\alpha$ and spinors $\phi$ on $X$,
defines a linear map, for which
we will use the shorthand notation
\[ V: (\phi,\alpha)
      \mapsto (\psi\circ\gamma)\cdot (\phi,\alpha)=(\phi^\tau, \alpha^\tau).
\]
In total we obtain bundle isomorphisms $\CA\to{\CA}^\tau$ and 
$\CC\to\CC^\tau$ of the Hilbert space bundles over the
identification $Pic^0(X)\stackrel{\cong}{\rightarrow} Pic^0( {X}^\tau)$ 
detailed above. All of these isomorphisms will be denoted by $V$.

\begin{thm} Gluing via the map $V$ induces for $b=b^+_2(X)$ an isomorphism
\[\pi^b_{{\T},H}(Pic^0(X);ind(D))\to 
\pi^b_{{\T},{H}^\tau}(Pic^0({X}^\tau);ind({D}^\tau)),\]
which identifies the classes of the monopole maps of $X$ and ${X}^\tau$
for corresponding $K$-theory orientations.
\end{thm}

\noindent
This theorem, the proof of which will be given in the next paragraph, implies
the gluing theorem stated in the introduction. This is a consequence of 
the next two propositions, applied to the case where at most one component of $X$
is not diffeomorphic to the standard four-sphere.

\begin{prop} The monopole map $\mu$ for a $K$-theory orientation on 
$X=\coprod_{i=1}^n X_i$ is the product of the monopole maps on the components of $X$
\[\mu=\prod_{i=1}^n\mu_i:\CA=\prod_{i=1}^n\CA_i\to\prod_{i=1}^n\CC_i=\CC.\]
Thus the associated stable equivariant cohomotopy element
is the smash product 
\[[\mu]=\wedge_{i=1}^n[\mu_i]\in\pi^{b^+_2(X)}_{{\T}^n,\, \oplus H_i}(Pic^0(X); ind(D))\]
of the cohomotopy elements associated to the respective components. 
The action of the torus ${\T}^n$ on the sum $\oplus_{i=1}^n H_i$ is factorwise. \qed
\end{prop}

This proposition is merely spelling out the obvious. The gluing
map $V$ is ${\T}$-equivariant with respect to the action of the diagonal subgroup of the torus. 

\begin{prop}
The stable cohomotopy element associated to any $\sc$-structure on a 
connected four dimensional
manifold $X$ with vanishing Betti numbers $b_1=b_2=0$ is the class of the 
identity map
\[[\mu]=[{\mathrm id}]\in \pi^0_{{\T}, H}(\ast)\cong {\Z}.\]
\end{prop}

\begin{pf}
In this case, the equivariant index ${\mathrm ind}(D)\in RO({\T})$ of the Dirac operator
is zero. Hence, the
ring $\pi^0_{{\T},H}(\ast)$ coincides with the Burnside ring $A({\T})\cong {\Z}$.
This isomorphism can be described as the map
\[
\pi^0_{{\T},H}(\ast)\to\pi^0_{st}(\ast)\cong{\Z}\]
induced by restriction to fixed point sets (cf. \cite{tD}, 133ff).
However, on the ${\T}$-fixed point set, the monopole map
is just the linear isomorphism
\[
{\mathrm d}+{\mathrm d}^\ast:\Omega^1(X)\to\Omega^2_+(X)\oplus\Omega^0(X)/{\R}.\]
\end{pf}


\section{Proof of the gluing theorem}

Let $\mu$ and ${\mu}^\tau$ denote the monopole maps on the $\sc$-manifolds
$X$ and ${X}^\tau$, respectively. The diagram
$$\begin{CD}   { \CA}     &   @>{\mu}>>   &    {\CC}          \\
       @V{V}VV&                 &    @VV{V}V      \\
   {\CA}^\tau& @>{{\mu}^\tau}>> & {\CC}^\tau
\end{CD}$$

\noindent
does not commute. The theorem claims that it commutes up to
suitable homotopy. What does suitable homotopy mean in this context?
The homotopy of course should be a homotopy through Fredholm maps, i.e. 
(nonlinear) compact deformations of linear Fredholm maps. However,
for the purposes in the present article these Fredholm maps
need not satisfy the boundedness condition of \cite{BauerFuruta} at all times. Instead, 
consider homotopies of Fredholm maps \[ {\mu}_t=l_t+c_t: {\CA}\to {\CC},\]
starting  from $\mu_0=\mu$ and ending at $\mu_1=V^{-1}{\mu}^\tau V$,
of the following kind:
There is a bounded disk bundle $D\subset {\CA}$ with bounding sphere
bundle $S$ over $Pic^0(X)$ such that both at the start and at the end of the homotopy
all solutions are contained in the disk bundle. During the homotopy no solution 
is allowed to cross the 
bounding sphere bundle. In more technical terms this 
means \[\mu_t^{-1}(0)\subset D\,\,\,\text{ for}\,\,\, t\in \{0,1\},\] and 
\[\mu_t^{-1}(0)\cap S = \emptyset\,\,\,\text{ for}\,\,\, t\in [0,1].\] 
When composing homotopies (as will be done below), 
the first condition, of course, need only be checked
at the very beginning of the first homotopy and at the end of the last homotopy.

\noindent
The reason why it suffices to consider this sort of homotopy was basically
mentioned in \cite{BauerFuruta}: 
The suitable homotopy is a homotopy of Fredholm maps
of pairs of spaces (over $Pic^0(X)$)
\[\mu_t: (D, S)\to ({\CC},{\CC}\setminus \{0\}).\]
Upon restriction to finite dimensional subbundles, the inclusions 
\[(D, S)\to ({\CA}^+, {\CA}^+\setminus D^\circ)\leftarrow 
({\CA}^+,\emptyset^+)\] and 
\[({\CC}^+,\emptyset^+)\to ({\CC}^+,({\CC}\setminus \{0\})^+)\leftarrow 
({\CC},{\CC}\setminus \{0\})\]
via homotopy equivalence and excision 
give isomorphisms of groups of pointed stable maps from
Thom spaces to spheres on the one side and stable maps  from pairs
(disk bundle, sphere bundle) to pairs (vector space, pointed vector space)
on the other side. Both groups describe the stable cohomotopy group
$\pi^{b^+_2(X)}_{{\T},H}(Pic^0(X);ind(D))$.

\noindent
In the proof the help of several homotopies in the above sense will be invoked:
The first two homotopies 
from $\mu$ (resp. ${\mu}^\tau$) to an auxiliary map 
$P:\CA\to\CC$ (resp. ${P}^\tau:{\CA}^\tau\to {\CC}^\tau$)
will tame the quadratic terms in the monopole map:
As an operator on sections over $X$, the map $P$
differs from the monopole map $\mu$ only over the
long neck. Over the short neck, the operators $P$ and $P^\tau$ 
are the linearisations of 
$\mu$ and $ \mu^\tau$.
The boundedness control during the homotopy is achieved by the use
of Weitzenb\"ock formulas for both the Dirac operator and the 
covariant derivative. Positivity of scalar and Ricci curvature,
respectively, on the long neck provide sufficient 
control on the spinor and
form components of solutions. The homotopy is set up so to keep  
the pivotal pointwise estimate on the spinor of a solution at all times.
The estimates on the forms will follow. However, in order to tune 
the estimates on spinors and forms,  
it may be necessary to stretch the long neck even longer. 
The final homotopy then ``rotates'' the operator $P$ into the operator $V^{-1}{P}^\tau V$.

\noindent
For $R\leq L$, let $\rho_R$ and $\rho_R^\tau$ be smooth cutoff 
functions, defined on $X$ and $X^\tau$ with values in $[0,1]$ of the following form:
The functions coincide, when restricted to the components $X_i^\pm$. On the middle
part $N(R-1)$ of length $2R-2$ of the long neck, the functions
vanish identically. Outside $N(R)$, the functions take the value $1$. In the remaining part 
$N(R)\setminus N(R-1)$, the function is constant on each sphere $\{r\}\times S^3$. 
The homotopies
\[\rho_{R,t}\stackrel{def}{=}(1-t)+t\rho_R,\]
for the time parameter $t$ in the unit interval, describe a homotopy from
the constant map $1$ to the function $\rho_R$ on $X$ and on ${X}^\tau$.

\subsection{The standard estimates}
The estimates in the proof of the main theorem are variations of the standard estimates 
used to show compactness of the moduli space. The argument
from \cite{KroMro} (compare \cite{BauerFuruta}, 3.1.), 
gives a norm estimate for  any point $(\phi, a)$
in the fibre over a point $A$ in the Picard torus  with $ \mu(A,\phi,a)   =  0  $.
Let's recall the three main steps of the argument:

{\bf Step 1:}   Applying the 
          Weitzenb\"ock formula for the Dirac operator  
           associated to the connection $A+a=A'$, one gets 
            a pointwise estimate:
            \[
                    \Delta|\phi|^2
                  \,\, \le \,\, 2<D^*_{A'}D_{A'}\phi-{s\over4}\phi+{1\over2}F_{A'}^+\phi,\phi>
              \]               
             The equalities $D_{A'}=0$ and $F_{A'}^+=\sigma(\phi)$ thus imply an estimate
            \[
                    \Delta|\phi|^2+{s\over2}|\phi|^2+{1\over2}|\phi|^4 \le 0
             \]       
        At the maximum of $|\f|^2$, its Laplacian is non-negative. So
             one obtains a pointwise estimate $ {s}|\f|^2+|\f|^4 \le 0
             $ for the norm of the spinor.

{\bf Step 2:} The Sobolev estimate $|a|_{C^0}\le C_s ||a||_{L^p_1}$ for some $p>4$ and
              the elliptic estimate 
              $||a||_{{\L}_1^p}\le C_e(||d^+a||_{{\L}_0^p}+||d^*a||_{{\L}_0^p}+||pr(a)||)$
              combine with the equality $d^+a=F^+_A+\sigma(\phi)$ to an estimate 

             \[
                |a|_{C^0}\le C_sC_e(||F^+_A||_{{\L}^p_0}+
                                     ||\sigma(\f)||_{{\L}^p_0}).
             \]

{\bf Step 3:}            
           For the bootstrapping assume for $i\le k$ inductively 
           ${\L}_{i-1}^{2p}$-bounds on $ ( \phi , a )$ with $p=2^{k-i}$.  
           To obtain ${\L}_{i}^p$-bounds, compute:
           \[   ||(\phi,a)||^p_{{\L}^p_i}  -   ||(\phi,a)||^p_{{\L}^p_0}
                             = 
                                 ||(D_A \phi,d^+a)||^p_{{\L}^p_{i-1}}                       
                            =      ||(a \phi,F_A^++\s(\f))||^p_{{\L}^p_{i-1}}.
           \]
           The summands in the last expression are bounded by the assumed   
           ${\L}_{i-1}^{2p}$-bounds on $ ( \phi , a )$.

\subsection{Varying the length of the long neck}
In the course of the proof it will be necessary to vary the 
length of the long neck. It is necessary to keep control on how 
the constants in step 2 of the argument do change.
To formulate the result needed, let $X'=X\setminus N(L-1)$ denote the 
complement of the middle part of length $2(L-1)$ of the neck. 

\begin{prop}\label{unabhaengig}
There is a constant $C_1$, independent of the length of the long neck, such that for any
smooth partition $\varphi_+, \varphi_-$ of unity on $X$, such that $\varphi_+ $ is identical $1$ on 
$X^+\cap X'$ and vanishes on $X^-\cap X'$, the following elliptic estimate holds:
$$|a|_{C^0(X')}\leq C_1(||(d^*+d^+)\varphi_+a||_{L^p}+||(d^*+d^+)\varphi_-a||_{L^p})+C||pr(a)||.$$
\end{prop}

\begin{pf}
The Sobolev constant $C_s$ depends only on the local geometry of the manifold: 
Choose for each $x\in X$ a bump function $\beta_x:X\to [0,1]$ with small support near $x$ 
and $\beta_x(x)=1$. The bump functions $\beta_x$ can be chosen such that 
their $C^1$-norm is bounded by a constant $M$, independently of the length
$2L$ of the neck. Let $c$ be a Sobolev constant which works for a particular $L$.  Then the inequality 
$$|a|_{C^0} = max_{x\in X}\,|\beta_x a|\le max_{x\in X}\,c||\beta_x a||_{{\L}^p_1}\le 
Mc||a||_{{\L}^p_1}$$ 
shows that the constant $C_s=Mc$ is independent of $L$.

It remains to gain control of the constants in the elliptic estimate. It is convenient
to invoke the theory of elliptic operators for manifolds with tubular ends, as
in \cite{Donaldson}. To do so, we attach
semi-infinite tubes to both positive and negative parts of $X=X^-\cup X^+$, to get manifolds 
$$
Y^-=X^-\cup ([0,\infty[\times S^3)\,\,\,\text{ and}\,\,\, Y^+=(]-\infty,0]\times S^3)\cup X^+
$$ 
with tubular ends. Let $f^-_{\alpha}\geq 0$ be a smooth function on $Y^-$ which 
vanishes identically on the complement of $N(L-2)$ and which grows with slope $\alpha$ 
towards the end, i.e. 
$f^-_\alpha(r,s)= \alpha (r+L-2) $ for $y=(r,s)\in {\R}\times S^3$ with $r\geq (-L+3)$.
A function $f^+_\alpha$ is defined similarly on $Y^+$ with $f^+_\alpha(r,s)= -\alpha(r-L+2)$ for
$r\leq (L-3)$. The weighted Sobolev spaces $L^{p,\alpha}_k(Y^\pm)$ are defined as the completions
of spaces of sections on bundles over $Y^\pm$ with respect to the norm 
$$
||g||_{L^{p,\alpha}_k(Y^\pm)}=||exp(f^\pm_\alpha)g||_{L^{p}_k(Y^\pm)}.
$$

After suitably identifying the bundles 
$$\Lambda^1(T^*({\R}\times S^3))
\cong 
{\R}\times(\Lambda^0\oplus\Lambda^1)T^*S^3
\cong
(\Lambda^0\oplus\Lambda^2_+)T^*({\R}\times S^3)$$
the operator $d^*+d^+$ gets the shape $\frac{\partial}{\partial t} +L$, where $L$ denotes the self adjoint
elliptic operator
$$L=
\left(
\begin{matrix}
0&d^*\\
d&*d
\end{matrix}
\right)
$$
on the three-sphere. As explained in \cite{Donaldson}, this operator $\frac{\partial}{\partial t} +L$
induces a Fredholm map $L^{p,\alpha}_1(Y^\pm)\to L^{p,\alpha}_0(Y^\pm)$ if $\alpha$ is not in the
spectrum of the operator $L$, whose square is the Laplace operator on the three-sphere. If $\alpha$
is negative and greater than the maximal negative eigenvalue of $L$, kernel and cokernel of the operator
$\frac{\partial}{\partial t} +L$ on $Y^\pm$ are isomorphic to $H^1(Y^\pm;{\R})$ and $H^0(Y^\pm;{\R})\oplus H^2_+(Y^\pm;{\R})$, respectively.
We may apply the estimate $$||b||_{L^{p,\alpha}_1(Y^\pm)}\leq C^\pm||(d^*+d^+)b||_{L^{p,\alpha}_0(Y^\pm)} + C||pr(b)||$$
to the compactly supported forms 
$b=\phi_\pm a$ on $Y^\pm$ which we get from a one-form $a$
an $X$. By construction both $a$ and $b$ have the same pointwise norm
along $X'$. The claim follows from this elliptic estimate and 
the Sobolev estimate with constant $C_1=C_s max(C^\pm)$.
\end{pf}

\subsection{The first homotopy}

\noindent
Consider the homotopy $\mu_t:\CA\to\CC$ defined by
\[\mu_t(A,\phi,a)=
(A,D_{A+a}\phi, F^+_{A+a}-\rho_{L,t}\sigma(\phi), pr(a), {\mathrm d}^\ast a).\]

\begin{lem}
The preimage $\mu_t^{-1}(0)$ is uniformly bounded for all times
$t\in[0,1]$.
\end{lem}

\begin{pf}
The standard argument applies with a minor change in step 1. At the maximum of 
$|\phi|^2$, the  estimate
$$ s|\phi|^2+\rho_{L,t}|\phi|^4\le  0$$
holds for the spinor component of a solution. The scalar curvature
is positive along the long neck. So the maximum is attained
outside the long neck and is bounded by the norm of the scalar curvature.
The norm of the term $\rho_{L,t}\sigma(\phi)$ in the rest of the argument
is bounded by a multiple of the norm of $\sigma(\phi)$. These bounds are 
independent of the parameter $t$.
\end{pf}

\subsection{The second homotopy}

\noindent
The next homotopy
moderates the second quadratic term over $N(2)$:
\[
\mu_{t+1}(A,\phi,a)= (A, D_{(A+\rho_{2,t}a)}\phi, F^+_{A+a}-\rho_L\sigma(\phi), pr(a), d^*a)
\]
This homotopy starts at $\mu_1$ and ends at $\mu_2=P$. Note that 
the latter differential operator is linear on the short neck.  
This second 
homotopy is more delicate than the first one. We need to
restrict it to a bounded disk as
explained in the beginning of this paragraph. 
In order to get the necessary
bounds on the solutions during the homotopy, it may, moreover, 
become necessary to 
stretch the long neck even longer, like playing the trombone. So we have to make sure that
the bounds we get for solutions $(\phi,a)$ to the operator $\mu_1$ are independent of the
length of the long neck. 
The spinor component causes no problems as it satisfies a $C^0$-bound
\[ |\phi|^2\leq S=\mathrm{max}(0,-s),\]
where $s$ denotes the scalar curvature of $X$ which is independent of $L$. 
The $C^0$-bound for the one-form
component of a solution used in step 2 of the standard argument was of the form
\[    |a|_{C^0}\leq  C_sC_e(||F^+_A||_{{\L}^p_0}+||\rho_L\sigma(\f)||_{{\L}^p_0})
               \leq  C_sC_e(|F_A^+|_{C^0}+S)vol(X')
\]
When stretching the neck, the factor $C_e$ is problematic, since it depends a priori on the
global geometry of $X$. So we will resort to make use of the less handy estimate \ref{unabhaengig}.

\begin{lem}\label{einsform}
There is a constant $U$ and a threshold length $L_0$ such that for any solution $(\phi,a)$ 
to the operator $\mu_1$ on the manifold $X(L)$ with neck of length $2L\geq 2L_0$ the one-form component 
satisfies the $C^0$-bound $|a|\leq U$.
\end{lem}

\begin{pf}
The one-form component of a solution on $X$ is harmonic on the part
$N(L-1)$ of the neck.
Because of nonnegative Ricci curvature along the neck, 
the maximum principle holds for the norm of such a one-form along $N(L-1)$.
Let $\varphi_+,\varphi_-$ be a partition of unity as in \ref{unabhaengig} with approximately
constant slope $|d\varphi_\pm|< L^{-1}$ along the neck. From \ref{unabhaengig} we get the estimate
\[|a|_{C^0(X')}\leq C_1(||(d^*+d^+)\varphi_+a||_{L^p}+||(d^*+d^+)\varphi_-a||_{L^p})+C||pr(a)||.\]
Since $pr(a)=0$, the latter summand does not contribute. We compute 
\[||(d^*+d^+)\varphi_\pm a||_{L^p}\leq ||\varphi_\pm(d^*+d^+)a||_{L^p}+||d\varphi_\pm||_{L^p}|a|_{C^0(N(L-1))}.\]
The vanishing of $(d^*+d^+)a$ along $N(L-1)$ and the maximum principle
\[|a|_{C^0(N(L-1))}\leq|a|_{C^0(X')}=|a|_{C^0(X)}\] thus lead to an estimate
\[|a|_{C^0(X)}\leq C_1(|F_A^+|_{C^0}+S)vol(X')+ C_1L^{\frac1p - 1}vol(S^3)|a|_{C^0(X)}.\]
The claim follows with $L_0=(2C_1vol(S^3))^{\frac43}$ and $U=2C_1(|F_A^+|_{C^0}+S)vol(X')$.
\end{pf}

\begin{lem}\label{Schwellenlemma}
If the long neck of $X$ is longer than the threshold length $L_1\geq L_0$, then the following holds:
Any solution  $(\phi, a)\in\mu^{-1}_t(0)$ for $1\leq t\leq 2$, which satisfies the 
$C^0$-bounds $|\phi|^2\leq 2S$ and $|a|\leq 2U$,  also satisfies the stricter 
$C^0$-bounds $|\phi|^2 \leq S$ and $|a|\leq U$. 
\end{lem}

\begin{pf}
Note that the one-form component of a solution on $X$ is again harmonic on the part
$N(L-1)$ of the neck and hence satisfies the maximum principle.
Moreover, because of the product structure of the neck,
such an harmonic one-form on the neck splits into a sum $a=a_i+a_s$ of harmonic one-forms,
according to the direct sum decomposition of the cotangent bundle.
The harmonic summand $a_s$ pointing in the sphere direction
satisfies an inequality
\[\Delta|a_s|^2\leq -2<Ric(a_s),a_s>\]
Since the Ricci tensor in direction of the sphere is positive definite,
\[-2<Ric(a_s), a_s>\,\,\leq \,-\delta^2 |a_s|^2\] for some $\delta>0$. 
If $\alpha$ denotes the sum $\sum_{i=1}^n|a_{s,i}|^2$ over all components of the
neck, then $\alpha$
satisfies a differential inequality
\[\frac{d^2\alpha}{dr^2}\ge-\Delta_X\alpha\ge \delta^2\alpha.\]
Thus $|a_s|^2$ is bounded by the function 
\[\frac{4nU^2\,\cosh (\delta r)}{\cosh(\delta(L-1))} .\] In particular, there is 
exponential decay of the norm of $a_s$ towards the middle of the neck.

Now we are ready to obtain for the spinor component of a solution a sharper bound than the assumed one.
The argument of step 1 works fine, if the spinor component attains its maximum outside $N(L-1)$.
So it suffices to make sure that the maximum is not attained along $N(L-1)$. Let's analyze, what may go wrong:
If the spinor component of a solution during the homotopy attains its maximum
in $N(L-1)$, then at that maximum, it satisfies an inequality
\[0\,\leq\, \Delta |\phi|^2 \,\leq\, -\frac{s}{2}|\phi|^2 + <(d\rho_{2,t}\smash a)^+\phi,\phi>.\]
Because of $d\rho_{2,t}\smash a = d\rho_{2,t}\smash a_s$, the norm of the latter summand
decays exponentially with the length of the long neck. If one stretches the long neck,
the second summand will decrease so that finally
the scalar curvature summand in the inequality will prevail. So if the neck is long enough,
the spinor component of a solution cannot attain its maximum in $N(L-1)$.

To get the sharper bound on the one-form component, one can now use exactly the same argument as in
\ref{einsform}, since the main ingredients for the argument hold: The bound on the spinor component is
unchanged, as is the fact that $a$ is harmonic along $N(L-1)$.
\end{pf}

\subsection{The third homotopy}
To finish the proof of the gluing theorem it remains to construct a homotopy 
between $P$ and $V^{-1}{P}^\tau V$. Note that both operators differ only over the
short neck in $X$. Because both differential operators are linear over the short neck,
their difference is a  multiplication operator:
\[V^{-1}{P}^\tau V=P+{\mathrm d}log(V)\]
For $t\in [0,1]$, consider the matrix valued function 
$\psi\circ t\gamma:[-L,L]\times S^3\to SO(n)$.
Multiplication of spinors or forms with this matrix valued function defines a map
$V_t$ over the long neck: Pairs of forms or spinors over the long 
neck are mapped to pairs of forms or spinors over the long neck. 
Multiplication with $\psi\circ t\gamma$ will not 
make sense outside the long neck. However, conjugation $V_t^{-1}PV_t=P+{\mathrm d}log(V_t)$
extends nicely to an operator $P+{\mathrm d}log(V_t):\CA\to \CC$ over all of $X$ for $0\leq t\leq 1$.
These operators
\[V_t^{-1}PV_t=P+{\mathrm d}log(V_t)\]
provide the final homotopy in the argument. 

\begin{lem} If the long neck of $X$ is longer than a threshold length 
$L_2\geq L_1$, then the following holds:
Any solution $(\phi,a)\in (P+{\mathrm d}log(V_t))^{-1}(0)$ for  $t\in[0,1]$, 
which is bounded by the  $C^0$-bounds $|\phi|^2\leq 2S$ and $|a|\leq 2U$, 
even satisfies the stricter $C^0$-estimates $|\phi|^2 \leq S$ and $|a|\leq U$. 
\end{lem}

\begin{pf}  
The proof makes use of the following observation: Let $(\phi,a)$ be a 
solution to the partial differential equation 
$$(P+{\mathrm d}log(V_t))(\phi,a)=V_t^{-1}PV_t(\phi,a)=0.$$
Then $(\phi,a)$ is also a solution to the operator $\mu_2$
over the complement of the short neck, since the $P+{\mathrm d}log(V_t)$ and $\mu_2$ only differ
over the short neck $N(1)$. Over the long neck, on the other hand,
the rotation $V_t^{-1}(\phi, \alpha)$ of $(\phi,a)$
is a solution to the operator $P$.

The proof of the $C^0$-bound for $|\phi|^2$ works almost as in \ref{Schwellenlemma}:
The Weitzenb\"ock argument can be applied
to $\phi$ on the complement of the short neck and to $V_t^{-1}\phi$ along $N(L-1)$
the same way, since $|\phi|^2=|V_t^{-1}\phi|^2$. Again we need to satisfy the condition that
the norm of $a_s$ along the short neck is small compared to the scalar curvature: 
Exponential decay of $|a_s|=|(V_t^{-1}a)_s|$ towards the middle of the neck 
follows as in \ref{Schwellenlemma} from the fact 
that $V_t^{-1}a$ as a solution to $P=\mu_2$ is harmonic along $N(L-1)$. Combined with
the assumed $C^0$-bound $|a|\leq 2U$, the condition again can be met by neckstretching.

Since $V_t^{-1}a$ is harmonic along $N(L-1)$, its $C^0$-norm $|V_t^{-1}a|=|a|$ obeys the maximum 
principle along this part of the neck. Let $\beta_\pm$ denote 
a smooth cutoff function supported on $X^\pm\setminus N(1)$
with $\beta_\pm(x)=1$ for $x$ in the complement of $N(L-1)\cap X^\pm$ and with
almost linear decay along the neck, i.e. $|d\beta_\pm|<\frac{2}{L}$. 
Then the function $\varphi=1-\beta_+-\beta_-$
is supported on $N(L-1)$. 
The elliptic estimate \ref{unabhaengig} gives a bound 
\begin{eqnarray*}
|a|_{C^0(X)}
&\leq& 
C_1(||(d^*+d^+)\beta_- a||_{L^p_0}+||(d^*+d^+)\varphi V_t^{-1}a||_{L^p_0}
+||(d^*+d^+)\beta_+ a||_{L^p_0})\\
&\leq& 
C_1(|F_A^+|_{C^0}+S)vol(X')+ C_1(||d\beta_-||_{L^p_0}+
||d\varphi||_{L^p_0}+||d\beta_+||_{L^p_0})|a|_{C^0(X')}\\
&\leq& 
\frac12 U+ 4C_1(\frac12L)^{\frac1p -1})vol(S^3)(2U).
\end{eqnarray*}
Again, the second summand can be made arbitrarily small by neckstretching.
\end{pf}

\noindent
For a proof of the theorem, one finally has to compose homotopies:
The first two homotopies combine to a
homotopy $\mu \sim P$. The third homotopy is between
$P$ and $V^{-1}P^\tau V$. Again the first two homotopies, conjugated by $V$, 
combine to the homotopy $V^{-1}\mu_{2-t}^\tau V$
from $V^{-1}P^\tau V$  to $V^{-1}\mu^\tau V$. The composition 
of these homotopies satisfies the suitability conditions harped on 
at the beginning of the section.


\section{Applications and Problems}
  
A $K$-theory orientation (or equivalently a $\sc$-structure)
on a four dimensional manifold is the same as a stably
almost complex structure. This is because
the natural map between the respective
classifying spaces $BU\to BSpin^c$,
has the appropriate connectivity. 
Such a stably almost complex structure comes with a
first Chern class, which is the integer lift $c\in H^2(X;{\Z})$ 
of the second Stiefel-Whitney class. 
The index of the Dirac operator associated to the $\sc$-structure
has complex dimension 
\[d={\mathrm ind}_{\C}(D)= \frac{c^2-sign(X)}{8}.\]
Via the Pontrijagin-Thom construction, the stable cohomotopy invariant
associated to a $K$-oriented four-manifold ideally can be thought of as encoding the  
equivariant framed bordism class of an ${\T}$-equivariant manifold $\tilde M$.
The quotient $M$ of this equivariant manifold by the ${\T}$-action then would be 
the moduli space considered in Seiberg-Witten theory. 
Actually, this ``ideal'' picture does not hold in general. The reason is that
the transversality arguments used in the Pontrijagin-Thom construction fail to
hold in an equivariant setting. As a consequence, the moduli spaces considered
in Seiberg-Witten theory do come with
singularities in those cases where equivariant transversality fails to hold.
Nevertheless, the moduli space will have an ``expected dimension''
\[k=2d-(b^+_2-b_1+1).\]
Let's consider an example where equivariant transversality fails. 
The following fact about equivariant maps is well known and can be proved
by the use of the equivariant $K$-theory mapping degree (cf. e.g. \cite{tD}):

\begin{prop}\label{Tammo}
Let $f:({\R}^n\oplus {\C}^{m+d})^+\to ({\R}^n\oplus {\C}^{m})^+$ be an ${\T}$-equivariant
map such that the restricted map on the fixed points has degree $1$. Then
$d\leq 0$ and $f$ is homotopic to the inclusion. \qed
\end{prop}
Such maps are not nullhomotopic, as their restriction to the fixed point set is not,
but the ``expected dimension'' $2d-1$ is negative. If equivariant transversality held
in this situation, then we get a contradiction: The manifold associated via 
Pontrijagin-Thom construction would be empty by dimension reasons. But the empty manifold
is associated to the nullhomotopic map.

The stable equivariant maps of \ref{Tammo} actually arise as monopole maps associated to 
four-manifolds with $b^+_2=b_1=0$, that is  
of manifolds with negative definite intersection form with vanishing first Betti number.
Applying the main theorem \ref{Summensatz} to this case results in a generalization of
the well-known ``blowing-up'' theorem:

\begin{cor} Let $X$ and $N$  be closed oriented four-manifolds with the 
intersection form on $N$ negative definite and $b_1(N)=0$. Fix $\sc$-structures
on both manifolds. The equivariant stable homotopy invariant of 
the connected sum is \[\mu_{X\#N}=\mu_X\wedge \gamma({\C})^{|d|}\]
with the diagonal ${\T}$-action, 
where $\gamma({\C})$ is the one-point compactified ${\T}$-map ${\C}^{0}\to {\C}^1$ and 
$d={\mathrm ind}_{\C}(D_N)\leq 0$. In particular, if $b^+_2(X)-b_1(X)>1$, then
the integer Seiberg-Witten invariants
of $X$ and $X\# N$ are the same, if $2|d|$ is not greater than the 
``expected dimension'' of the monopole moduli space associated to $X$.
\end{cor}
The condition $b^+_2(X)-b_1(X)>1$ can actually be removed by considering maps ``up to equivariant
homotopy modulo the fixed point set'' as explained in \cite{BauerFuruta}.

\begin{pf} Because of \ref{Tammo}, the first part is just a restatement of the gluing
theorem in this special case. The last statement is trivial in case the expected dimension
$k$ of the monopole moduli space associated to $X$ is odd, since then the 
Seiberg-Witten invariants of both $X$ and $X\#N$ vanish by definition. Otherwise the claim
follows by the very definition of the comparison map \cite{BauerFuruta}, 3.3.
from the stable cohomotopy invariants to the integers: Let 
$\mu_X\in \pi^b_{{\T}}(Pic^0(X);ind(D))$ denote the stable cohomotopy invariant of $X$.
Then $\mu_X\wedge \gamma({\C})^{\frac{k}2}$ as an element in 
$\pi^b_{{\T}}(Pic^0(X);ind(D)-\frac{k}2{\C})\cong {\Z}$ 
is just the Seiberg-Witten invariant of $X$. 
\end{pf}

By a theorem of Hirzebruch and Hopf, the expected dimension is zero if and only if
the stably almost complex structure on $X$ is actually an almost complex structure.
Indeed, in all currently known (at least to the author)
examples of four dimensional $\sc$-manifolds 
with nonvanishing integer valued
Seiberg-Witten invariant, the $\sc$-structure 
is associated to an almost complex structure.

\begin{notitle}\label{Taubes}Here are some general facts about Seiberg-Witten invariants: 
\begin{itemize}
\item
By a theorem of Taubes \cite{Taubes}, the Seiberg-Witten invariant of a symplectic 
four-manifold is, up to sign convention, 1.
\item
If $X$ is a four-dimensional manifold underlying a K\"ahler surface, then solutions
to the monopole equations correspond to holomorphic sections of certain complex line
bundles, compare \cite{FriMor}. As a consequence, the stable cohomotopy invariants and
a fortiori the Seiberg-Witten invariants, are nonzero at most for such $K$-orientations
which correspond to almost complex structures. In particular, the moduli space either has 
expected dimension $0$ or is empty.
\end{itemize}
\end{notitle}

For simplicity, consider from now on only manifolds with vanishing first Betti number.
If the $K$-orientation of $X$ is associated to an almost complex structure, then
the stable map $\mu$ nonequivariantly is
an element of the stable homotopy group $\pi_1^{st}(S^0)$.
This group has two elements, the trivial map and the 
Hopf map $\eta$. 

\begin{prop}\label{Hopf}The cohomotopy invariant of an almost complex
manifold with vanishing first Betti number 
nonequivariantly is the Hopf map if and only if 
both $b^+_2$ is congruent $3\, mod\, 4$ and the integer valued
Seiberg-Witten invariant is odd.
\end{prop}

\begin{pf} As was shown in \cite{BauerFuruta}, the stable
        cohomotopy invariant 
for an almost complex manifold $X$ is
$SW(X)\kappa_d$, where
$\kappa_d$ is a generator of the stable cohomotopy group
$\pi^{2d-2}({\C}P^{d-1})\cong{\Z}\kappa_d$ and $2d=b^+_2+1$. 
The composition 
\[S^{2d-1}\to S^{2d-1}/{\T}={\C}P^{d-1}\buildrel{\kappa_d}\over\rightarrow S^{2d-2},\]
of $\kappa_d$ with the quotient map of the free ${\T}$-action
is the nontrivial Hopf element $\eta\in\pi_1^{st}(S^0)$ iff
$d$ is even. This can be seen as follows: The quotient map
is the attaching map of the top cell in ${\C}P^d$. If $z$
denotes the generator in $H^2({\C}P^d;{\F}_2)$, 
then the Steenrod square $Sq^2(z^{d-1})$
is nonzero if and only if $d$ is even. But $Sq^2$ detects
the Hopf map. 
\end{pf}

By considering products of the ${\T}$-equivariant Hopf map, we get the following result, which
in particular implies \ref{power4}

\begin{prop}\label{Produkte}
Let $X$ be a connected sum of $n\geq 2$ almost 
complex manifolds $X_i$ with vanishing first Betti numbers. 
The stable equivariant cohomotopy element of $X$ is nonvanishing
if and only if the following conditions are satisfied:
\begin{itemize}
\item
For each summand, $b^+_2(X_i)$ is congruent $3\, mod\, 4$. 
\item
The integer valued Seiberg-Witten invariants are odd for each summand.
\item
If $n\geq 4$, then $n=4$ and $b^+_2(X)$ is congruent $4\, mod\, 8$.
\end{itemize}
\end{prop}

\begin{pf} 
The result follows from the Atiyah-Hirzebruch spectral sequence, which computes 
the groups $\pi^{2d-1-n}({\C}P^{d-1})$, which contain the stable cohomotopy 
invariants in these cases. The results of this computation are stated in
\cite{BauerFuruta}, 3.5.\\
First consider the case $n=2$.
The ``if'' part of the statement follows from the fact that
the square of the Hopf map is the only 
nonzero element in the second stable stem. \\
On the other hand, the group $\pi^{2d-3}({\C}P^{d-1})$
is nonzero if and only if $d$ is even. In this case
its only nontrivial element is the composition of the quotient map 
${\C}P^{d-1}\to {\C}P^{d-1}/{\C}P^{d-2}=S^{2d-2}$
with the Hopf map. Combined with \ref{Hopf}, the result follows.
\\
The result for $n=3$ is immediate from the result for $n=2$ 
and from the fact that nonequivariantly the cube of the Hopf map is nonzero.
\\
If $n=5$, the stable cohomotopy invariant is identified with a torsion element
in the group $\pi^{2d-6}({\C}P^{d-1})$, which is torsion free. So in particular,
the stable cohomotopy invariant is zero, if $n\geq 5$.
\\
It remains to consider the case $n=4$. This is more delicate, because the 
fourth power of the Hopf map is zero nonequivariantly. As is well known,
the product in stable homotopy theory can be defined equivalently in two 
ways. One way is by smash product, the other is by composition.
If we use the second description, we see that the first Hopf map factors
through  the quotient map 
$S^{2d-1}\to {\C}P^{d-1}$ by the ${\T}$-action, followed
by the projection to the top cell. It remains to consider, whether projection to
the top cell, followed by the cube of the nonequivariant Hopf map is stably nonzero
as a map from ${\C}P^{d-1}$ to $S^{2d-5}$. But this can be read off the Atiyah-Hirzebruch
spectral sequence: This map is nonzero if and only if the torsion subgroup in 
$\pi^{2d-5}({\C}P^{d-1})$ is divisible by $8$, which is the case if and only if
$d$ is divisible by $8$, which is the case if and only if  $b^+_2(X)$ is congruent $4\, mod\, 8$.
\end{pf}

\begin{pf}(of \ref{Separatisten2})
The proof follows from \ref{Produkte} and \ref{Taubes}:
The stable cohomotopy invariant of the connected sum
$X$ of two symplectic manifolds with $b^+_2\equiv 3\,mod\, 4$
nonequivariantly lies in the second stable stem
$\pi_2^{st}(S^0)$
and actually is the
square of the Hopf map. If $X$ splits into a connected sum
$X\cong X_1\#X_2$ of manifolds with $b^+_2\equiv 1\,mod\,2$, then
the  stable cohomotopy invariants of the respective summands
are elements of the $n_1$-th and $n_2$-th stable stem with
$n_i$ both odd and $n_1+n_2=2$. Since their product is
the nonvanishing element in the second stable stem, both 
$n_i$ are equal $1$. By the Hirzebruch-Hopf theorem 
both manifolds are almost complex.
Now \ref{Produkte} applies again to yield the contradiction.
\end{pf}

\begin{pf}(of \ref{Separatisten3})
If $X$ is the connected sum of three symplectic manifolds with
 $b^+_2\equiv 3\,mod\, 4$, then the stable cohomotopy invariant
is the cube of the Hopf map in the third stable stem
(nonequivariantly). If $X$ splits into a connected sum 
$X\cong X_1\#X_2$ with $b^+_2(X_1)\equiv 1\,mod\,2$, then (because
 of \ref{Hopf}) the manifold $X_1$ cannot be almost
 complex. Its stable cohomotopy invariant lies in an odd
 stable stem, but not in the first. The stable cohomotopy invariant of $X_2$ lies
 in an even stable stem. As their product is the nontrivial
 cube of the Hopf map in
 the third stable stem, only one possibility remains: The
 stable cohomotopy invariant of $X_2$ is an odd element in
 the zero'th stable stem $\pi_0^{st}(S^0)\cong \Z$. Now we
 have to apply equivariant homotopy theory \cite{tD} again: 
The stable cohomotopy invariant of $X_2$ is
 represented by an ${\T}$-equivariant map 
 \[f:({\R}^n\oplus {\C}^{m+\frac{b}2})^+\to
 ({\R}^{n+b}\oplus {\C}^{m})^+.\]
If $b=b^+_2(X_2)$ is nonzero, then the nonequivariant degree of such a
map is necessarily zero. This finishes the proof.
\end{pf}

\begin{pf}(of \ref{K3})
For the K3-surface the SW-invariants are completely known: They vanish
except for the one $\sc$-structure which lifts to a $Spin$-structure. For 
this the value is $1$, up to sign convention. From the statements above
it follows that for a connected sum $K\# K\# X$ or $K\#X$ for a 
simply connected K\"ahler surface $X$, the $\sc$-structures supporting
nontrivial 
stable cohomotopy invariants of the connected sum correspond to exactly
those $\sc$-structures on $X$ having odd SW-invariants (if $b^+_2$ is of
correct modulus). 
\end{pf}

\noindent
For simply connected manifolds, the $\sc$-structures are detected
by the first Chern classes of the associated spinor bundles. 
 The argument above can be rephrased the following way: For a 
connected sum of K\"ahler surfaces, the nontrivial stable cohomotopy
invariants detect the pairs (or triples) consisting of cohomology
classes of the summands which support odd Seiberg-Witten invariants.
This can be applied in special situations to recognize the summands:

\begin{pf}(of \ref{elliptic}) The proof is based on the known classification
of elliptic surfaces (see e.g. \cite{FriMor}) via SW-invariants.
To a simply connected minimal 
elliptic surface with given $b^+_2$ one can associate 
a pair $1\leq m\leq n$ of coprime integers which, 
together with the geometric 
genus $p_g=\frac{b^+_2-1}{2}$, classify the diffeomorphism type.
Note that $b^+_2$ is congruent to $3$ mod $4$ iff the geometric genus is
odd. 
The point in the proof is that one can recognize $m,n$ and $p_g$ for 
odd geometric genus from the pattern of the cohomology classes
corresponding to odd SW-invariants, the ``recognizable'' classes.
Here comes a description, how this can be accomplished combinatorially.

\noindent
The cohomology classes corresponding to nontrivial
SW-invariants are multiples of an indivisible element $f$ in the second
cohomology  with integer values of the elliptic surface. 
The multiplicities are of the form
$(p_g-1-2a)mn+(m-2b-1)n+(n-2c-1)m$ for nonnegative integers 
$a< p_g,
\,\,b< m$
and $c< n$. The value of the SW-invariant for such a multiple of $f$ 
is $ \left(\begin{matrix}
p_g-1\\a\end{matrix}\right)$. 

\noindent
Note that the distribution of basic classes is symmetric around the 
origin and the SW-invariant for the largest such multiple $k$ of $f$, 
where $a=b=c=0$, is odd. 
So there are at least two recognizable classes except in the case
of a K3-surface $p_g=m=n=1$, where there is
exactly one recognizable class.

\noindent
If there are no more than three recognizable 
classes, then either $m=n=1$ or  $n=2m=2p_g=2$.
In the latter case, the largest multiple is 1, in the former, it is
$p_g-1$, which is even.

\noindent
In the case of at least four recognizable classes, consider the second
but largest multiple. In case $n>1$, the integer $m$, which is
half the difference, is coprime to the largest multiple. If $m=n=1$,
then there is an $0<2a\leq p_g-1$ with $ \left(\begin{matrix}
p_g-1\\a\end{matrix}\right)$ odd. This integer $a$ has to be even,
because otherwise  $ \left(\begin{matrix}
p_g-1\\a+1\end{matrix}\right)$, which is obtained from it by multiplication
with a rational number having $a+1$ in its denominator, could not be
an integer. In this case, half the difference cannot be coprime to
the largest multiple; both are even. 
This makes it possible to distinguish the
$m=n=1$-cases.

\noindent
Finally consider the largest multiple $k$ and the second largest multiple.
It can be assumed that half the difference is coprime to $k$
and thus equals $m$. Consider the multiples $(k-2\lam m)f$ for $\lam \ge 1$.
The SW-invariants associated to these classes will be 1 for 
$\lam < n$ and zero or $p_g-1$, anyway even, for $\lam = n$. 
This characterizes the second
integer $n$. Knowing both integers this way, the geometric genus
follows from the formula for $k$.

\noindent
In the situation of a connected sum of no more than three elliptic surfaces,
as in \ref{elliptic}, the number of summands can be read off the dimension
of the moduli spaces having nontrivial invariants. The cohomology
classes associated to nontrivial invariants are situated in a 
bounded region in a sublattice of the second cohomology of rank at most 3.
They form a box and the pattern characterizing the individual 
summands can be found on the respective edges of the box.
\end{pf}

\begin{notitle}
The monopole map $\mu$ may be pertubed quite a bit without changing
the resulting stable cohomotopy invariant. 
For example the term $\sigma(\phi)$ may be
replaced by some function $ f(|\phi|)\sigma(\phi)$ as long as $f$
does not decay too fast at infinity. Any polynomial 
with positive leading coefficient will give different moduli
spaces, but the same stable cohomotopy invariant.\end{notitle}

\noindent
Let's end this article with two problems, which grow out of considering
stable cohomotopy invariants of connected sums. The first problem is related
to a well known theorem of C.T.C. Wall stating that 
any two homeomorphic simply connected four dimensional differentiable 
manifolds become diffeomorphic after taking connected sum with finitely
many copies of $S^2\times S^2$. Moreover, in many cases of algebraic
surfaces it is known that it suffices to take connected sum with only
one such copy of $S^2\times S^2$.

\begin{problem} Suppose $X$ and $Y$ are homeomorphic, simply connected
differentiable four-manifolds. Do they become diffeomorphic after
taking connected sum with sufficiently many K3-surfaces? 
\end{problem}

\begin{problem} Are there manifolds realizing stable cohomotopy
elements in $\pi^i({\C}P^{d-1})$ other than the powers of
the Hopf map $\eta$, for example the element 
associated to the stable Hopf map $\nu:S^7\to S^4$? 
Or more generally: Is there an indecomposable $K$-oriented 
four-manifold with nonvanishing stable cohomotopy invariant,
which is not almost complex?
\end{problem}


\end{document}